\documentclass [11pt] {article}
\usepackage{amsfonts}
\usepackage{amssymb}
\usepackage{times}
\usepackage{graphicx}

\newtheorem{lemma}{Lemma}[section]
\newtheorem{corollary}{Corollary}[section]
\newtheorem{theorem}{Theorem}

\newtheorem{definition}{Definition}[section]
\def\limo{\lim_{\omega}}

\def\<{\langle}
\def\>{\rangle}

\def\bH{{\bf H}}

\def\bS{{\bf S}}

\def\to{\rightarrow}

\def\bX{{\bf X}}
\def\bY{{\bf Y}}
\def\bA{{\bf A}}
\def\bH{{\bf H}}

\def\bD{{\bf D}}
\def\bm{{\bf \mu}}
\def\rk{{\rm rk}}
\begin{document}

\title{Higher order fourier analysis as an algebraic theory I.}
\author{{\sc Bal\'azs Szegedy}}

\maketitle

\abstract{Ergodic theory, Higher order fourier analysis and the hyper graph regularity method are three possible approaches to Szemer\'edi type theorems in abelian groups. In this paper we develop an algebraic theory that creates a connection between these approaches. Our main method is to take the ultra product of abelian groups and to develop a precise algebraic theory of higher order characters on it. These results then can be turned back into approximative statements about finite Abelian groups.
}

\tableofcontents

\section{Introduction}

This paper is an announcement of a new approach to higher order Fourier analysis. It contains proofs and some explanations but a longer, more detailed version is on the way. The next paper will contain further results on the structure of higher order characters and will have a special emphasis on the finite translations of our results. Note that every result on ultra products in this paper automatically translates to an approximative finite result with epsilons.

In a paper \cite{ESz} by G. Elek and myself we develop a measure theoretic approach to dense hypergraphs. The central theorem of the paper creates a correspondence between growing sequences of hypergraphs and various measurable sets. A powerful feature of that approach is that extra algebraic structures (such as groups structures) on the hypergraphs can be detected on the limit objects and thus can be related to the regularity lemma. The general reason is that ultra products preserve axiomatizable structures.

Our starting point is the following question. Let $S$ be a subset of a finite abelian group $A$ and let $H$ be the $k$ uniform hypergraph consisting of the $k$-edges
$$\{(a_1,a_2,\dots,a_k)~|~a_1+a_2+\dots+a_k\in S\}.$$
{\it What is the regularization of $H$ in terms of the structure of $S$?}

\noindent More generally, if $f:A\rightarrow\mathbb{C}$ is a bounded function, how can we regularize $f(x_1+x_2+\dots+x_k)$?

In the paper \cite{ESz} we prove a certain uniqueness of regularization. This suggest that there must be a nice algebraic answer to the above problem.
It turns out that for $k=2$ ordinary Fourier analysis of $f$ gives a full answer to the question (when dominant terms in the Fourier expansion are considered).
However for $k>2$ a higher order version of Fourier analysis is needed.
Such a theory was started first by Gowers in \cite{Gow2} and later continued by Green, Host, Kra, Tao, Ziegler and many others (see reference list).

Our contribution to the subject is a clean algebraic version of Higher order Fourier analysis (on the ultra product of finite abelian groups) which is quite similar to ordinary Fourier analysis. For every $n$ there are order $n$ characters on our ultra product Abelian group that form an orthonormal basis for certain Hilbert spaces depending on $n$. These characters can be described as generators of rank-one modules. Further more the $n$-th order characters are forming an Abelian group with respect to point wise multiplication that we call the $n$-th dual group.

We emphasize that on a finite Abelian group ordinary linear characters form a basis of the dual space and this phenomenon leaves no room for higher order characters that are fully orthogonal to the ordinary ones. However on the ultra product group the first order characters are not forming a full basis. This allows us to develop our algebraic theory.

Note that our results can be extended to compact Abelian groups but we will do it in a forthcoming paper. Also the relationship between higher order characters and nil-systems will be discussed in a the second part of this paper.

Since translating the results back to finite statements is a routine thing we don't focus on it. We give the details in a longer version of the same paper.

Finally we mention that our language allows us to define a version of higher order representation theory for non-commutative groups but we don't know applications yet.

\subsection{Explanation of the Results and further remarks}

Is this section we give a short explanation of the results in the paper.

Let $\{A_i\}_{i=1}^\infty$ be an increasing sequence of finite Abelian groups and let $\bA$ be the ultra product of them. Using the language from \cite{ESz} there is a $\sigma$-algebra $\mathcal{A}$ and a shift invariant probability measure $\bm$ on $\mathcal{A}$.

Let $\lambda_i:A_i\rightarrow\mathbb{C}$ a sequence of linear characters. It is easy to see that the ultra limit of these functions is a linear character on $\bA$. Functions arising this way bill be called {\bf first order characters} on $\bA$.
It turns out that first order characters are not forming a basis in $L_2(\bA,\mathcal{A},\bm)$. There is a smaller sigma algebra $\mathcal{F}_1$ such that they are an orthonormal basis in $L_2(\bA,\mathcal{A},\bm)$. We call $\mathcal{F}_1$ the {\bf first order Fourier $\sigma$-algebra}.

A result in this paper gives a simple characterization for the elements in $\mathcal{F}_1$. A measurable set $M$ is in $\mathcal{F}_1$ if and only if the $\sigma$-algebra generated by the shifts $M+x$ of $M$ is separable. In other words countable many shifts of $M$ generate all the shifts. We call such sets $M$ {\bf separable elements} in $\mathcal{A}$. Note that $\mathcal{F}_1$ is not a separable $\sigma$-algebra despite of the fact that every element of it is separable.

The {\bf second order Fourier analysis} comes from an interesting phenomenon. The are measurable sets $M$ in $\mathcal{A}$ that are not separable but countable many shifts of $M$ together with $\mathcal{F}_1$ generate all the shifts of $M$. Such elements are called {\bf relative separable} elements with respect to $\mathcal{F}_1$. These elements are forming a $\sigma$-algebra which is denoted by $\mathcal{F}_2$.

Continuing this process $\mathcal{F}_n$ denotes the $\sigma$-algebra of relative separable elements with respect to $\mathcal{F}_{n-1}$. We say that $\mathcal{F}_n$ is the $n$-th order Fourier $\sigma$-algebra.

\medskip

Our definition of {\bf Higher order Forier analysis} is the following: {\it Decompose $L_2(\mathcal{F}_n)$ into irreducible $L_\infty(\mathcal{F}_{n-1})$ modules!}

\medskip

The {\bf main theorem} in this paper says that $L_2(\mathcal{F}_n)$ can uniquely be decomposed into {\bf rank one $L_\infty(\mathcal{F}_{n-1})$ modules} that are pairwise orthogonal to each other.
Every such module is generated by a function $\lambda$ which takes complex values of absolute value one. These functions are measurable in $\mathcal{F}_n$ and furthermore $f(x):=\lambda(x)\overline{\lambda(x+t)}$ is measurable in $\mathcal{F}_{n-1}$.
The function $\lambda$ us unique up to multiplication with elements from $L_\infty(\mathcal{F}_{n-1})$. Note that $\mathcal{F}_0$ is defined to be the trivial $\sigma$-algebra. This shows that the first order Fourier analysis gives back the classical theory.

Our language is basically {\bf ergodic theoretic}. One can look at $\bA$ as a measure preserving system acting on itself. Since the $\sigma$-algebras $\mathcal{F}_i$ are shift invariant, they can be regarded as {\bf factors} of the system. Our main theorem says in this formulation that $\mathcal{F}_n$ is an Abelian extension of $\mathcal{F}_{n-1}$.
A lots of intuition for this paper (especially the part about cubic structures) was borrowed from a beautiful paper by Host and Kra \cite{HKr}.

In a continuation of this paper we will analyze these systems more with a special emphasis on their relationship with nil-systems.

We mention that our theory can be looked at from the {\bf hypergraph regularization} point of view: If $f:\bA\rightarrow\mathbb{C}$ is a bounded function then the ``first level'' regularization of the $n$-variable function $f(x_1,x_2,\dots,x_n)$ is obtained by taking the projection $g=E(f|\mathcal{F}_{n-1})$ and then considering $g(x_1+x_2+\dots+x_n)$.

The so called {\bf Gowers's norm} $U_n(f)$ is $0$ if and only if $f$ is orthogonal to $\mathcal{F}_{n-1}$. In other words for every bounded function $f:\bA\rightarrow\mathbb{C}$ there is a (unique!) decomposition
$$f=g+f_1+f_2+f_3+\dots$$
converging in $L_2$ such that $U_n(g)=0$ and the functions $f_i$ are contained in different rank-one modules over $L_\infty(\mathcal{F}_{n-2})$.

We mention that a part of this theory can be developed for non-commutative groups. One can take their ultra product and then consider that tower of relative separable $\sigma$-algebra.
It turns out that they are {\bf isometric extensions} of each other.
Details will be worked out later.

\subsection{Characteristic Factors}

In this short chapter we explain an interesting phenomenon. The concept of a {\bf characteristic factor} is crucial in all the three theories: Ergodic theory, Higher order Fourier analysis and Hypergraph regularization. Later on, when we relate the three theories it will be apparent.

Characteristic factors are classical in Ergodic theory. It turns out that certain averages in measure preserving systems remain the same when the set gets projected to a courser $\sigma$-algebra. The smallest such algebra is characteristic for the type of average.

In hypergraph theory they show up in the following way. The so called counting lemma in the infinite setting \cite{ESz} says that if we want to count sub-hypergraphs in a $k$-uniform hypergraph then it is enough to project the the indicator function of the edge set to certain courser $\sigma$-algebras. The coursest such $\sigma$-algebra gives rise to hypergraph regularization.

As it is pointed out in the present paper, Higher order Fourier analysis can also be interpreted in a similar way. There we project functions to the Fourier $\sigma$-algebras.

\section{Manipulating with $\sigma$-algebras}

\subsection{Basics and notation}

Throughout this paper we fix a non-principal ultra-filter $\omega$ on the natural numbers $\mathbb{N}$. For a sequence of sets $\{X_i\}_{i=1}^\infty$ we denote their ultra product by $[\{X_i\}_{i=1}^\infty]$. We will also use the shorthand notation $[\{X_i\}]$ once the running index $i$ is clear from the context. For a bounded sequence $a_1,a_2,\dots$ of complex numbers we denote their ultra limit by $\limo a_i$.

\subsection{Measure and $\sigma$-topology on ultra-products}

Let $\omega$ be a fixed ultrafilter on $\mathbb{N}$. Let $\{X_i\}_{i=1}^\infty$ be an infinite sequence of finite sets. We denote their ultra product $[\{X_i\}_{i=1}^\infty]$ by $\bX$.

\begin{definition} A subset $\bH\subseteq\bX$ is said to be {\bf perfect} if
$$\bH=[\{H_i\}_{i=1}^\infty]$$
for some subsets $H_i\subseteq X_i$.
A subset $\bH\subseteq\bX$ is said to be {\bf open} if it is the union of at most countable many perfect sets.
\end{definition}

We denote by $\mathcal{P}(\bX)$ the collection of perfect sets and by $\mathcal{O}(\bX)$ the open sets on $\bX$. We say that the open sets  define a $\sigma$-{\bf topology} on $\bX$. For $\bH\in\mathcal{P}(\bX)$ we introduce the measure $\bm(\bX)$ as
$$\bm(\bH)=\lim_\omega(|H_i|/|X_i|).$$
The sigma algebra generated by $\mathcal{P}(\bX)$ on $\bX$ will be denoted by $\mathcal{A}=\mathcal{A}(\bX)$. The measure $\mu$ extends to a $\sigma$-additive probability measure on $\mathcal{A}$.

\begin{lemma}[countable compactness]\label{sigcom} If $\bA$ is the union of countable many open sets $\{\bH_i\}_{i=1}^\infty$ then there is a natural number $n$ such that $\bA=\cup_{i=1}^n\bH_i$.
\end{lemma}

\begin{definition} A function $f:\bX\rightarrow\mathbb{C}$ is called {\bf perfect} if
$$f([{x_i}_{i=1}^\infty])=\limo f_i(x_i)$$
for some sequence of functions $f_i:X_i\rightarrow\mathbb{C}$.
\end{definition}

\subsection{Perfect $\sigma$-algebras}

Let $\{X_i\}_{i=1}^\infty$ be a sequence of finite sets and let $\{2^{X_i}\}_{i=1}^\infty$ be the sequence of their power sets.
The elements of the ultra product $$[\{2^{X_i}\}_{i=1}^\infty]$$ are in a one to one correspondence with $\mathcal{P}(\bX)$ so by abusing the notation we can identify the two objects.
The advantage of this representation of $\mathcal{P}(\bX)$ is that it creates a $\sigma$-topology on $\mathcal{P}(\bX)$.

\begin{definition} The sigma algebra $\mathcal{B}\subseteq\mathcal{A}(\bX)$ is called {\bf perfect} if it is generated by a set in $\mathcal{P}(\mathcal{P}(\bX))$.
\end{definition}

An important advantage of perfect $\sigma$-algebras is that the conditional expectation operator behaves in a nice measurable way.

\begin{definition} Let $\mathcal{B}\subseteq\mathcal{A}$ be a $\sigma$-algebra. We say that $\mathcal{B}$ has the {\bf projection property} if there is a $\mathcal{A}(\mathcal{P}(\bX)\times \bX)$-measurable function
$$f:\mathcal{P}(\bX)\times \bX\rightarrow [0,1]$$
such that for every $\bH\in\mathcal{P}(\bX)$
$$f(\bH,x)=E(\bH|\mathcal{B})(x).$$
\end{definition}

\begin{lemma} Any perfect $\sigma$-algebra has the projection property.
\end{lemma}

\begin{proof} Assume that $\mathcal{B}$ is generated by the ultra product of the sequence $\{S_i\}_{i=1}^\infty$ where $S_i\subseteq 2^{X_i}$. Let $\bS$ be the ultra product of the sequence $\{S_i\}_{i=1}^\infty$. For every natural number $n$ and $i$ we construct functions $f_{n,i}:2^{X_i}\times X_i\rightarrow [0,1]$ in the following way.
Let us fix a subset $H\subseteq X_i$ and take the projections of $H$ to all finite $\sigma$-algebras generated by at most $n$ elements from $S_i$. Each such projection is a function of the form $X_i\rightarrow [0,1]$.
We define the function $x\mapsto f_{n,i}(H,x)$ as one of the projections that has minimal $L_2$ norm among all possible projections.
We denote by $\hat{f}_n$ the ultra limit function of the sequence $\{f_{n,i}\}_{i=1}^\infty$.
One can observe that for every fixed perfect set $\bH$ the function $x\rightarrow\hat{f}_n(\bH,x)$
is a projection of $\bH$ to a finite $\sigma$-algebra generated by $n$-sets from $\bS$. Furthermore its norm is minimal among all such projections.
This implies that the measurable function $\lim_{n\to\infty}\hat{f}_n$ has the desired property.
\end{proof}

\subsection{$\sigma$-algebras in hyper-graph regularization and Gowers's norms}

We review some of the language developed in \cite{ESz} to regularize hypergraphs with ultra product.
Let $k$ be a fixed natural number and $S$ be a subset of $\{1,2,3,\dots,k\}$. We denote by $P_S$ the natural projection from $\bX^k$ to $\bX^S$.
The pre-image of $\mathcal{A}(\bX^S)$ under $\bX^k$ will be denoted by $\sigma(S)$. Furthermore $\sigma(S)^*$ denotes the $\sigma$-algebra generated by all the $\sigma$-algebras $\sigma(S')$ where $S'$ is a proper non-empty subset of $S$.

A set $\bH\in\mathcal{A}(\bX^k)$ is called {\bf quasi-random} if the projection
$E(\bH|\sigma([k])^*)$ is the constant function.
Let $f:\bX^k\rightarrow\mathbb{C}$ be an $L_2$ function. The function $f$ is called quasi-random
if $E(f|\sigma([k]^*)$ is the zero function.
The nice thing about quasirandom functions is that they are forming a Hilbertspace.
Gowers ... introduced the so called {\bf octahedral-norms} the characterize quasirandom finctions in the finite setting. A similar lemma applies in the infinite setting but the epsilon's are disappears.

\begin{definition} The octahedral norm $O_k(f)$ is the $2^k$-th root of the integral
$$\int_{x_{i,j}\in\bX}\prod_{c\in \{0,1\}^k}f(x_{1,c_1},x_{2,c_2},\dots,x_{k,c_k})^{q(c)}$$
where $1\leq\i\leq k$ ~,~ $j\in\{0,1\}$ and $q(c):\mathbb{C}\rightarrow\mathbb{C}$ denotes the conjugation operator raised to the power $\sum_{i=1}^k c_i$.
\end{definition}

\begin{lemma} The function $f:\bX\rightarrow\mathbb{C}$ is quasi-random if and only if $O_k(f)=0$.
\end{lemma}

\begin{definition}[Gowers's norms]\label{gowers} Let $f:\bX\rightarrow\mathbb{C}$ be an $L_2$-function and let $f_k:\bX^k\rightarrow\mathbb{C}$ denote the function
$f_k(x_1,x_2,\dots,x_k)=f(x_1+x_2+\dots+x_k)$. Then the Gowers norm $U_k(f)$ is defined as $O_k(f_k)$.
\end{definition}

\subsection{Weak Orthogonality and coset $\sigma$-algebras}

Let $H\subseteq L_2(\bX,\mu)$ be a Hilbert space of functions on $\bX$. We denote by $P_H$ the orthogonal projection to $H$. We will need the following lemma

\begin{lemma}[Integration in Hilbert spaces]\label{inth} Let $f:\bX\times\bY\rightarrow\mathbb{C}$ be an $\mathcal{A}(\bX\times\bY)$ measurable function such that $g_y(x)=f(x,y)$ is in a Hilbert space $H$ for every $y\in\bY$. Then
$$g(x)=\int_y f(x,y)~d\bm_2$$
is in $H$.
\end{lemma}

\begin{proof} Following the next equations
$$\|g(x)\|^2=\int_x\int_y f(x,y)\overline{g(x)}=\int_y\int_x f(x,y)\overline{g(x)}=$$
$$=\int_y\int_x f(x,y)\overline{P_H(g)(x)}=\int_y\int_x f(x,y)\overline{P_H(g)(x)}=$$
$$=\int_x g(x)\overline{P_H(g)(x)}=\|P_H(g)\|^2$$

we get that $g=P_H(g)$ which shows that $g\in H$.
\end{proof}

\begin{definition} Two Hilbert spaces $H_1$ and $H_2$ in a bigger Hilbert space $H$ will be called {\bf weakly orthogonal} if for any $f\in H_1$ we have that $P_{H_2}(f)\in H_1$. Let $\mathcal{B}_1$ and $\mathcal{B}_2$ be two sub $\sigma$-algebras in the same measure space. We say that $\mathcal{B}_1$ is weakly orthogonal to $\mathcal{B}_2$ if for any $f\in\mathcal{B}_1$ the function $E(f|\mathcal{B}_2)$ is measurable in $\mathcal{B}_1$.
\end{definition}

\begin{lemma} Weak orthogonality is a symmetric relation.
\end{lemma}

\begin{proof}
Weak orthogonality for Hilbert spaces is clearly a symmetric notion. It
means that the orthogonal complements of $H_1\cap H_2$ in $H_1$ and $H_2$ are orthogonal. Weak orthogonality of the $\sigma$-algebras $\mathcal{B}_1$ and $\mathcal{B}_2$ is equivalent with the weak orthogonality of $L_2(\mathcal{B}_1)$ and $L_2(\mathcal{B}_2)$ so it is again symmetric.
\end{proof}

\begin{definition} Let $\{A_i\}_{i=1}^\infty$ be an infinite sequence of finite Abelian groups. We denote their ultra product $[\{A_i\}_{i=1}^\infty]$ by $\bA$. A subgroup $\bH\subseteq\bA$ is called {\bf perfect} if there are subgroups $H_i\subseteq A_i$ with $\bH=[\{H_i\}_{i=1}^\infty]$. The {\bf coset $\sigma$-algebra} corresponding to a perfect subgroup $\bH$ is the subalgebra $\mathcal{A}(\bA,\bH)\subseteq\bA$ consisting of the measurable sets which are unions of cosets of $\bH$. A $\sigma$-algebra $\mathcal{B}\subseteq\mathcal{A}$ is called $\bA$-{\bf invariant} (or invariant) if for every $\bY\in\mathcal{B}$ and $a\in\bA$ we have that $\bY+a\in\mathcal{B}$.
\end{definition}

Clearly, any coset $\sigma$-algebra is invariant.
The advantage of coset $\sigma$-algebras is that the projection operator to them can be computed by an averaging operator.
Let $f:\bA\rightarrow\mathbb{C}$ be an $L_2$ function. Then we can define its average $T(f)$ according to $\bH$ by
$$T(f)(x)=\int_{h\in\bH}f(x+h).$$
Lemma \ref{inth} implies that if $f$ is measurable in an invariant $\sigma$-algebra $\mathcal{B}$ then $T(f)$ is also measurable in $\mathcal{B}$.
It is clear from the definition that $(T(f),w)=(f,w)$ for every $w\in L_2(\mathcal{A}(\bA,\bH))$ and that $T(f)\in L_2(\mathcal{A}(\bA,\bH))$. We obtain that $T(f)$ coincides with the projection
$E(f|\mathcal{A}(\bA,\bH))$. As a consequence we get the next lemma.

\begin{lemma}\label{weako} A coset $\sigma$-algebra on $\bA$ is weakly orthogonal to any invariant $\sigma$-algebra.
\end{lemma}

\subsection{Modules and relative separability}

\begin{definition} Let $\mathcal{B}\subseteq\mathcal{A}(\bX)$ be a $\sigma$-algebra and let $R=L_{\infty}(\mathcal{B})$ be the ring of bounded $\mathcal{B}$-measurable functions $f:\bX\mapsto\mathbb{C}$ with the point wise multiplication. A Hilbert space $H\subseteq L_2(\bX)$ is an $R$-{\bf module} if $hr\in H$ for every $h\in H$ and $r\in R$. The {\bf rank} ${\rm rk}(H)$ is the minimal size of a subset $S\subseteq H$ such that the $R$-module generated by $S$ is dense in $H$.
\end{definition}

\begin{lemma}\label{modvet} Let $H_1$ and $H_2$ be two $R$ modules. Then $P_{H_2}(vr)=rP_{H_2}(v)$ whenever $v\in H_1$ and $r\in R$.
\end{lemma}

\begin{proof} We have that
$$(P_{H_2}(vr),w)=(vr,w)=(v,w\bar{r})=(P_{H_2}(v),w\bar{r})=(rP_{H_2}(v),w)$$ holds for every $w\in H_2$.
By setting $w=P_{H_2}(vr)-rP_{H_2}(v)$ we get that $(w,w)=0$ which completes the proof.
\end{proof}

\begin{lemma}\label{modrk} Let $H_1\subseteq H_2$ be two $R$ modules. Then $\rk(H_1)\leq \rk(H_2)$.
\end{lemma}

\begin{proof} Lemma \ref{modvet} implies that if $S$ is a generating set of $H_2$ as an $R$ module then $S'=\{P_{H_2}(s)~|~s\in S\}$ is a generator set of $H_1$.
\end{proof}

\begin{lemma} Let $H_1\subseteq H_2$ be $R$-modules and let $H_3$ be the orthogonal complement of $H_1$ in $H_2$. Then $H_3$ is an $R$-module and
$E(f\bar{g}|\mathcal{B})$ is the $0$ function for every $f\in H_1$ and $g\in H_3$.
\end{lemma}

\begin{proof} We have that if $r\in R$, $h\in H_3$ and $t\in H_1$ then $(t,rh)=(t\bar{r},h)=0$ and so $rh\in H_3$. Let $e=E(fg|\mathcal{B})$. We have that lemma \ref{modvet} that $E(f\bar{g}\bar{e}|\mathcal{B})(x)=|e(x)|^2$ and using $(f,ge)=0$ we obtain that $\int|e(x)|^2=0$.
This shows that $e$ is the constant $0$ function.
\end{proof}

\begin{definition} Let $\mathcal{B}_1,\mathcal{B}_2\subseteq\mathcal{A}(\bX)$ be $\sigma$-algebras. We say that $\mathcal{B}_2$ is {\bf relative separable} over $\mathcal{B}_1$ if $\langle\mathcal{B}_1,\mathcal{B}_2\rangle$ can be generated by $\mathcal{B}_1$ and at most countable many extra elements as a $\sigma$-algebra.
\end{definition}

If $\mathcal{B}_1\subseteq\mathcal{B}_2$ then relative separability means that $L_2(\mathcal{B}_2)$ is an at most countable rank $L_\infty(\mathcal{B}_1)$-module.

\begin{lemma}[Independent complement] Let $\mathcal{B}_1\subseteq\mathcal{B}_2\subseteq\mathcal{A}(\bX)$ be two $\sigma$-algebras such that $\mathcal{B}_2$ is relative separable and relative atom less over $\mathcal{B}_1$. Then there is a separable $\sigma$-algebra $\mathcal{B}_3\subseteq\mathcal{B}_2$ which is independent from $\mathcal{B}_1$ and $\langle\mathcal{B}_1,\mathcal{B}_3\rangle=\mathcal{B}_2$.
\end{lemma}

\begin{lemma}[Relative basis] Let $\mathcal{B}_1\subseteq\mathcal{B}_2\subseteq\mathcal{A}(\bX)$ be two $\sigma$-algebras such that $\mathcal{B}_2$ is relative separable and relative atom-less over $\mathcal{B}_1$.
Then there is a system of at most countably many functions $f_i:\bX\rightarrow\mathbb{C}$ such that
$$E(f_i\overline{f_j}|\mathcal{B}_1)(x)=\delta_{i,j}.$$
\end{lemma}

\begin{lemma} Let $\mathcal{A}(\bA,\bH)$ be a coset $\sigma$-algebra, $\mathcal{B}$ be a sub $\sigma$-algebra in it and $\mathcal{B}_2\subseteq\mathcal{A}(\bA)$ a $\sigma$-algebra such that $\mathcal{B}$ is relative separable over $\mathcal{B}_2$. Then $\mathcal{B}$ is relative separable over $\mathcal{A}(\bA,\bH)\cap\mathcal{B}_2$.
\end{lemma}

\subsection{Slices of measurable sets}

Let $\bH\subseteq\bX^k$ be a subset and let $f:\bX^k\rightarrow\{0,1\}$ be its characteristic function. In this section we study the ``slices'' of $\bH$ according to the last coordinate. Let $x$ be an element in $\bX$. The $x$-slice of $\bH$ is the set
$$\bH_x:=\{(x_1,x_2,\dots,x_{k-1}~|~(x_1,x_2,\dots,x_{k-1},x)\in \bH\}.$$
The $\sigma$-algebra generated by all $x$-slices will be denoted by $\mathcal{S}(\bH)$.

\begin{lemma}\label{slice} The set $\bH$ is in $\sigma([k])^*$ if and only if $\mathcal{S}(\bH)$ is relative separable over $\sigma([k-1])^*$.
\end{lemma}

\begin{proof} First we assume that $\bH\in\sigma([k])^*$. Then there are countable many sets in $A_{i,j}\in\sigma([k]/\{j\})$ with $i\in\mathbb{N}$ and $1\leq j\leq k$ such that $\bH$ is measurable in the $\sigma$-algebra generated by $\{A_{i,j}\}$. If $j<k$ then $\mathcal{S}(A_{i,j})\in\sigma([k-1]^*)$ furthermore $\mathcal{S}(A_{i,k})$ is generated by one element. This shows that
$$\mathcal{S}(\bH)\subseteq\langle\mathcal{S}(A_{i,j})~|~i\in\mathbb{N},1\leq j\leq k\rangle$$
is relative separable over $\sigma([k-1])^*$.

For the other direction let $\mathcal{B}$ be a sigma algebra containing $\mathcal{S}(\bH)$ such that $\mathcal{B}$ is relative separable and relative atom less over $\sigma([k-1])^*$. Such a $\sigma$-algebra can be constructed by extending $\mathcal{S}(\bH)$ by at most countable many generators. Let $g_1,g_2,\dots$ be a relative basis of $\mathcal{B}$ over $\sigma([k-1])^*$.
We have that
$$f(x_1,x_2,\dots,x_k)=\sum_i g_i(z)E(f_x(z)\overline{g_i(z)}|\sigma([k-1])^*)$$
where $z=(x_1,x_2,\dots,x_{k-1})$.
This formula show that $f$ is measurable in $\sigma([k])^*$.
\end{proof}

Note that $\mathcal{S}(\bH)$ might depend on a $0$ measure change of $\bH$. For this reason we introduce a similar but more complicated notion. Let us denote the $\sigma$-algebra
$$\langle \sigma(S)~|~S\subseteq[k],~k\in S,~|S|=k-1\rangle$$ by $\beta([k])$.

\begin{definition} The {\bf essential} $\sigma$-algebra $\mathcal{E}(\bH)$ is the smallest $\sigma$-algebra containing $\sigma([k-1])^*$ in which all the functions
$$g(x_1,x_2,\dots,x_{k-1})=\int_{x_k}f(x_1,x_2,\dots,x_k)t(x_1,x_2,\dots,x_k)$$
are measurable where $t$ is an arbitrary bounded function in
$$L_\infty(\beta([k])).$$
\end{definition}

In the next lemma $\psi$ denotes the projection prom $\bX^k$ to $\bX^{k-1}$ which cancels the last coordinates.

\begin{lemma}\label{esprop} The $\sigma$-algebra $\mathcal{E}(\bH)$ has the following properties
\begin{enumerate}
\item $\mathcal{E}(\bH)=\mathcal{E}(\bH_2)$ whenever $\bH\triangle\bH_2$ has measure $0$
\item $\mathcal{E}(\bH)$ is relative separable over $\sigma([k-1])^*$
\item $\bH$ is almost measurable in $\langle \psi^{-1}(\mathcal{E}(\bH)),\beta([k])\rangle$.
\item the functions $f_{x_k}:=f(x_1,x_2,\dots,x_k)$ are measurable in $\mathcal{E}(\bH)$ for almost all $x_k\in\bX$.
\end{enumerate}
\end{lemma}

We will need the following lemma.

\begin{lemma}\label{hosszu} Let $\mathcal{B}_1,\mathcal{B}_2\subseteq\mathcal{A}$ be two weakly orthogonal $\sigma$-algebras on a probability space $(X,\mathcal{A},\mu)$ and $W$ be a subset measurable in $\langle\mathcal{B}_1,\mathcal{B}_2\rangle$.
Let furthermore $\mathcal{B}_3$ be the $\sigma$-algebra generated by all the functions
$E(ft|\mathcal{B}_1)$ where $f$ is the indicator function of $W$ and $t$ is in $L_\infty(X,\mathcal{B}_2)$. Then $H$ is almost measurable in $\langle\mathcal{B}_3,\mathcal{B}_2\rangle$.
\end{lemma}

\begin{proof}
First we prove the statement in the case when $\mathcal{B}_1$ and $\mathcal{B}_2$ are independent.
Let $W_2$ be a set with $|W_2\triangle W|\leq\epsilon$ of the form
$$W_2=\bigcup_{i=1}^n (A_i\cap B_i)$$
where $\{B_i\}_{i=1}^n$ is a partition of $X$ into $\mathcal{B}_2$ measurable sets and $A_i\in\mathcal{B}_1$ for every $1\leq i\leq n$.
The existence of such a $W_2$ is guaranteed by the fact that $W$ is measurable in $\langle\mathcal{B}_1,\mathcal{B}_2\rangle$.
We have that
$$\chi(W_2)=\sum_{i=1}^n \frac{1}{|B_i|}\chi(B_i)E(\chi(W_2)\chi(B_i)|\mathcal{B}_1).$$
where $\chi$ denotes the indicator function of a set.
Let $$g=\sum_{i=1}^n \frac{1}{|B_i|}\chi(B_i)E(\chi(W)\chi(B_i)|\mathcal{B}_1).$$
Now $$\|\chi(W_2)-g\|_1\leq\bigl\|\sum_{i=1}^n \frac{1}{|B_i|}\chi(B_i)E((|\chi(W_2)-\chi(W)|)\chi(B_i)|\mathcal{B}_1)\bigr\|_1.$$
Therefore by the independence of $\mathcal{B}_1$ and $\mathcal{B}_2$ we obtain that $\|\chi(W_2)-g\|_1\leq\epsilon$. Together with $\|\chi(W_2)-f\|_1\leq\epsilon$ we obtain that
$\|g-f\|_1\leq 2\epsilon$. Since $g$ is measurable in $\langle\mathcal{B}_3,\mathcal{B}_2\rangle$ the result follows.

Now we go to the general case. According to a lemma in... we can construct a $\sigma$-algebra $\mathcal{A}\subseteq\mathcal{B}_2$ such that $\mathcal{A}$ is independent from $\mathcal{B}_1$ and
$W$ is also measurable in $\langle\mathcal{B}_1,\mathcal{A}\rangle$. This completes the proof.
\end{proof}

Now we are ready to prove lemma \ref{esprop}

\begin{proof}
The first property is obvious from the definition.
From lemma \ref{inth} we obtain that $\mathcal{E}(\bH)\subseteq\mathcal{B}_2$. Then lemma \ref{modrk} shows the second condition.

For proving the third statement we use the fact that $\sigma([k-1])$ and $\beta([k])$ are weakly orthogonal. Then lemma \ref{hosszu} implies the statement.
\end{proof}

\subsection{Relative separable elements and Fourier $\sigma$-algebras}

Let $\{A_i\}_{i=1}^\infty$ be an infinite sequence of finite Abelian groups. We denote their ultra product $[\{A_i\}_{i=1}^\infty]$ by $\bA$. Lemma \ref{sigcom} implies that $\bA$ is countable compact in the $\sigma$-topology $\mathcal{O}(\bA)$.
Let $s:\bA\times\bA\rightarrow\bA$ denote the addition map defined by $s(a,b)=a+b$. In general it is not true that $s$ is continuous in the product $\sigma$-topology $\mathcal{O}(\bA)\times\mathcal{O}(\bA)$ on $\bA\times\bA$ but it will be obviously continuous in the $\sigma$-topology $\mathcal{O}(\bA\times\bA)$.

It is clear that the $\sigma$-algebra $\mathcal{A}$ and the measure $\bm$ is invariant under the action of $\bA$

\begin{definition} Let $\mathcal{B}$ be an arbitrary sub $\sigma$-algebra of $\mathcal{A}$ and let $\bH$ be an element of $\mathcal{A}$. We say that $\bH$ is {\bf relatively separable} with respect to $\mathcal{B}$ if there are countable many translates $T=\{\bH+a_i~|~i\in\mathbb{N},a_i\in\bA\}$ of $\bH$ such that the $\sigma$-algebra generated by $T$ and $\mathcal{B}$ is dense in $\sigma$-algebra generated by all the translates $\{\bH+a~|~a\in\bA\}$ and $\mathcal{B}$.
\end{definition}

The collection of all relative separable elements of $\mathcal{B}$ is again a $\sigma$-algebra which contains $\mathcal{B}$. We denote it by $\Upsilon(\mathcal{B})$.
It is clear that if $\mathcal{B}$ is invariant then so is $\Upsilon(\mathcal{B})$.
This means that iterating the operator $\Upsilon$ starting with the trivial $\sigma$-algebra we get an interesting sequence of invariant algebras.

\begin{definition}[Fourier $\sigma$-algebras] Let $\mathcal{F}_0$ denote the trivial $\sigma$-algebra on $\bA$ consisting of the empty set and the whole set $\bA$. We define the sequence of $\sigma$-algebras $\mathcal{F}_i$ recursively in the way that $\mathcal{F}_i=\Upsilon(\mathcal{F}_{i-1})$.
\end{definition}

Clearly all the $\sigma$-algebras $\mathcal{F}_i$ are invariant under the action of $\bA$.

\subsection{Characterization of the Fourier $\sigma$-algebras}

In this chapter we give equivalent characterizations of the $\sigma$-algebras $\mathcal{F}_k$.
Let $\bD_k$ denote the subgroup in $\bA^k$ consisting of elements $(x_1,x_2,\dots,x_k)$ with $x_1+x_2+\dots+x_k=0$. The factor group $\bA^k/\bD_k$ is isomorphic to $\bA$ where the isomorphism is given by $\tau_k(x_1,x_2,\dots,x_k)=x_1+x_2+\dots+x_k$. It is also clear that $\tau_k$ creates a measure preserving equivalence between the $\sigma$-algebras $\mathcal{A}(\bA)$ and $\mathcal{A}(\bA^k,\bD_k)$. We prove the following fact which is crucial in this theory.

\begin{lemma}\label{measeq} The map $\tau_k$ gives a measure preserving equivalence between $\mathcal{A}(\bA^k,\bD_k)\cap\sigma([k])^*$ and $\mathcal{F}_{k-1}$.
\end{lemma}

The proof will require the next simple lemma.

\begin{lemma}\label{seged} Let $\bY$ be a perfect subgroup in a Abelian group $\bA$. Assume that a measurable set $\bH\in\mathcal{A}(\bA,\bY)$ is relative separable over an invariant $\sigma$-algebra $\mathcal{B}\subseteq\mathcal{A}(\bA)$.
Then $\bH$ is also relative separable over $\mathcal{B}\cap\mathcal{A}(\bA,\bY)$.
\end{lemma}

\begin{proof} Let $T=\{\bH+t_i\}_{i=1}^\infty$ be a system of countable many translates of $\bH$ such that all other translates are in the $\sigma$-algebra generated by $T$ and $\mathcal{B}$. Let $a\in\bA$ be an arbitrary element and $\epsilon>0$ be a positive number. Then there is a finite expression of the form
$$\bH_\epsilon=\cup_{i=1}^n (A_i\cap B_i)$$
with $\bm(\bH\triangle\bH_{\epsilon})\leq\epsilon$
such that $A_i$ is measurable in the $\sigma$-algebra generated by $T$ and $B_i$ is in $\mathcal{B}$ for every $1\leq i\leq n$. Now we have that
$$g:=E(\bH_\epsilon|\mathcal{A}(\bA,\bY))=\sum_{i=1}^n\chi(A_i)E(B_i|\mathcal{A}(\bA,\bY)).$$
Furthermore $\|g-\chi(\bH)\|_2\leq\epsilon$ sice the distance of two functions is decreased by projections. From lemma \ref{weako} it follows that the functions $E(B_i|\mathcal{A}(\bA,\bY))$ are all measurable in $\mathcal{B}\cap\mathcal{A}(\bA,\bY)$. This completes the proof.
\end{proof}

Now we are ready to prove lemma \ref{measeq}.

\begin{proof} We prove the statement by induction on $k$. The case $k=1$ is obvious from the definitions. We assume that it is true for $k-1$. Let $\bH$ be a set in $\mathcal{A}(\bA)$ and $f$ its characteristic function. We have to prove that $\bH$ is measurable in $\mathcal{F}_{k-1}$ if and only if $f_k$ is measurable in $\sigma([k-1])^*$. (Recall that $f_k$ is the characteristic function of $\tau_k^{-1}(\bH)$.)

First assume that $\bH$ is measurable in $\mathcal{F}_{k-1}$. This means that $\bH$ is relative separable over $\mathcal{F}_{k-2}$. The slice $\sigma$-algebra $\mathcal{S}(f_k)$ consists of the translates of $f_{k-1}$ and so it is relative separable over $\tau_{k-1}^{-1}(\mathcal{F}_{k-2}$. By the induction hypothesis $f_{k-1}$ is also relative separable over $\sigma([k-1])^*$. This implies by lemma \ref{slice} that $f_k$ is measurable in $\sigma([k])^*$.

For the other direction we assume that $f_k$ is in $\sigma([k])^*$. It is clear that $\mathcal{E}(f_k)$ is an invariant $\sigma$ algebra on the Abelian group $\bA^{k-1}$. Since almost all the slices of $f_k$ are measurable in $\mathcal{E}(f_k)$ and furthermore all slices are translates of $f_{k-1}$ we get from the invariance of $\mathcal{E}(f_k)$ that all the translates of $f_{k-1}$ are measurable in $\mathcal{E}(f_k)$.
By Lemma \ref{esprop} we obtain that $f_{k-1}$ is a relative separable element over $\sigma([k-1])^*$.
It follows from lemma \ref{seged} that $f_{k-1}$ is relative separable over $\mathcal{A}(\bA^{k-1},\bD_{k-1})\cap\sigma([k-1])^*$. Applying $\tau_{k-1}^{-1}$ and the induction hypothesis, $\bH$ is relative separable over $\mathcal{F}_{k-2}$.
\end{proof}

As a corollary we get the next theorem.

\begin{theorem}\label{fouequiv} Let $f:\bX\rightarrow\mathbb{C}$ a function in $L_2$ and let $f_k$ be as in definition \ref{gowers}. Then the following statements are equivalent.
\begin{enumerate}
\item $f$ is measurable in $\mathcal{F}_{k-1}$
\item $f$ is orthogonal to every $g$ with $U_k(g)=0$.
\item $f_k$ is measurable in $\sigma([k])^*$
\item $f_k$ is orthogonal to every $g$ with $O_k(g)=0$.
\end{enumerate}
\end{theorem}

\begin{lemma}\label{projprop} The $\sigma$-algebras $\mathcal{F}_i$ have the projection property for every $i$.
\end{lemma}

\subsection{Group extensions}

Let $G$ be a second countable compact topological group with the Borel $\sigma$-algebra $\mathcal{G}$ on it. We denote the Haar measure on $G$ by $\nu$.
A function $f:\bA\rightarrow G$ is called a {\bf pre-cocycle} of order $k$ if $f$ is measurable in $\mathcal{F}_k$ and $g_t(a):=f(a+t)f(a)^{-1}$ is measurable in $\mathcal{F}_{k-1}$ for every fixed $t\in\bA$.

Two pre-cocycles $f$ and $f_2$ of order $k$ are called equivalent if $f_2=f_3fz$ where $z$ is a fixed element in $G$ and $f_3$ is measurable in $\mathcal{F}_{k-1}$.

Let us define the action of $\bA$ on $\bA\times G$ by
$$a.(b,g)=(b+a,f(b+a)f(b)^{-1}g)$$ and the action of $G$ on $\bA\times G$ by
$$h.(b,g)=(b,gh).$$
The two actions commute with each other.
It is clear that the above actions preserves the measure $\mu\times\nu$ where $\mu$ is the ultra product measure on $\bA$.
Furthermore the $\sigma$-algebra $\mathcal{F}_{k-1}\times\mathcal{G}$ is invariant under $\bA$.

Let $\mathcal{I}$ denote the $\sigma$-algebra of $\bA$ invariant sets in $\mathcal{F}_{k-1}\times\mathcal{G}$.
The action of $G$ leaves $\mathcal{I}$ invariant and acts on it in an ergodic way.
Now Lemma 3.23 in [GLASNER] implies that there is a closed subgroup $H$ in $G$ such that the action of $G$ on $\mathcal{I}$ is isomorphic to the action of $G$ on the left coset space of $H$ in $G$.
We say that $H$ is the Mackey group corresponding to $f$.  Note that it is defined up to conjugacy.

We say that $f$ is {\bf minimal} if $\mathcal{I}$ is trivial.
We introduce the averaging operator
$$Q:L_2(\mathcal{F}_{k-1}\times\mathcal{G})\rightarrow L_2(\mathcal{F}_{k-1}\times\mathcal{G})$$ by
$$(Qf)(b,g)=\int_a a.(b,g)~d\mu.$$
Lemma \ref{inth} shows that $Q$ has the required form.
It is also clear that $Q$ is the orthogonal projection to $L_2(\mathcal{I})$.

\begin{lemma}\label{findmin} There is a minimal pre-cocycle $f_2:\bA\rightarrow H$ such that $f$ and $f_2$ are equivalent where $H$ is the Mackey group of $f$.
\end{lemma}

\begin{proof} See [GLASNER] 3.25 pg 73
\end{proof}

\begin{lemma}\label{cosigma} If $f:\bA\rightarrow G$ is a minimal pre-cocycle then there is a $\sigma$-algebra $\mathcal{F}_{k-1}\subset\mathcal{H}\subset\mathcal{F}_k$ and a measure algebra equivalence
$\phi:\mathcal{F}_{k-1}\times\mathcal{G}\rightarrow\bA$ which commutes with the $\bA$ action.
Furthermore $\phi$ restricted to $\mathcal{F}_{k-1}\times Triv(G)$ gives an equivalence between $\mathcal{F}_{k-1}$ and $\mathcal{F}_{k-1}\times Triv(G)$ where $Triv(G)$ denotes the trivial $\sigma$ algebra on $G$.
\end{lemma}

\begin{proof} For a Borel set $B\subseteq G$ let $T(B)$ denote the set $\{g~|~f(g)\in B\}$.
We show that $\mu(M\cap T(B))=\mu(M)\nu(B)$ for any two measurable sets $M\in\mathcal{F}_{k-1}$ and $B\in\mathcal{G}$. Since $f$ is minimal we have that $Q(H\times B_2)$ is (almost everywhere) the constant function with value $\mu(M)\nu(B_2)$ for any Borel set $B_2$. Lets choose a countable dense set $\hat{B}$ of open sets in $\mathcal{G}$ with respect to the Haar measure. For almost every pair $(a,g)\in\bA\times G$ we have that
$$\int_b \chi_{M\times B_2}(b.(a,g))=\mu(M)\nu(B_2)$$
for every $B_2\in \hat{B}$ where $\chi_{M\times B_2}$ is the characteristic function of $M\times B_2$. Let $(a,g)$ be one such pair. then $-a.(a.g)=(0,f(-a)f(a)^{-1}g)=(0,g_2)$ is clearly another such pair. The value of the above integral for $(0,g_2)$ is equal to $\mu(M)\nu(B_2)$. On the other hand it is the same as $\mu(M\cap T(B_2g_2^{-1}f(0)))$. Using the translation invariance of $\nu$ this formula implies that $\mu(M)\nu(B_2g_2^{-1}f(0))=\mu(M\cap T(B_2g_2^{-1}f(0)))$. Since $\hat{B}$ is dense we get the formula for every $B\in\mathcal{G}$.

Now it is easy to check that $\phi(M\times G)=M\cap T(B)$ extends to a measure preserving equivalence and it is given by $$\phi(W)=\{a~|~(a,f(a))\in W\}.$$
\end{proof}

\begin{definition} Let $\mathcal{F}_{k-1}\mathcal{H}\mathcal{F}_k$ be a $\bA$ invariant $\sigma$-algebra. Then an automorphism of the measure algebra $(\mathcal{H},\mu)$ is called a {\bf k- automorphism} if $\sigma$ commutes with the action of $\bA$ on itself and leaves $\mathcal{F}_{k-1}$ invariant. The set of $k$-automorphisms will be denoted by $Aut_k(\mathcal{H})$. Any $k$ automorphism induces an action on $L_2(\mathcal{H})$ such that $E(\sigma(f)|\mathcal{F}_{k-1})=E(f|\mathcal{F}_{k-1})$.
\end{definition}

\begin{lemma}\label{action} Let $f$ be a minimal pre-cocycle of order $k$ and let $\mathcal{H}$ be the $\sigma$-algebra constructed in lemma \ref{cosigma}. Then $G$ induces a faithful action $\varrho:G\rightarrow Aut_k(\mathcal{H})$.
\end{lemma}

\begin{proof} The proof is an immediate corollary of lemma \ref{cosigma}. The action of $G$ on $\bA\times G$ induces a similar action on $\mathcal{H}$.
\end{proof}

\subsection{Cubic structure}

In this section we borrow part of our language from \cite{HKr}.
For a natural number $k$ let $V_k$ denote the set of all subsets in $[k]$. We will think of the elements of $V_k$ as the vertices of a $k$-dimensional cube.
We will use the notion of a $d$-dimensional ${\bf face}$ in the natural way.
Let us define the subgroup $B_k$ in $\bA^{V_k}$ as the collection of vectors $(a_i)_{i\in V_k}$ satisfying all the equations
$$a_p-a_q+a_r-a_s=0$$
where $\{p,q\}$ is a $1$ dimensional face and $\{p,q,r,s\}$ is a $2$ dimensional face of $V_k$.
We call the vertex corresponding to the empty set the $0$-vertex. If $v\in V_k$ is a vertex then the {\bf spider} $S(v)$ is the subset consisting of $v$ and its neighbors in $V_k$.
The following lemma shows that the projection of $B_k$ to the coordinates in $S(v)$ gives a bijection between $B_k$ and $\bA^{S(v)}$. In particular $B_k$ is isomorphic to $\bA^{k+1}$.

\begin{lemma}\label{spidrep} Let $\delta_k:\mathcal{A}^{S(0)}\mapsto\mathcal{A}^{V_k}$ denote the homomorphism defined by
$$\delta_k(a_0,a_1,a_2,\dots,a_k)_v=a_0+\sum_{i\in v} a_i.$$
Then $\delta_k$ gives a measure preserving isomorphism between $\bA^{k+1}$ and $B_k$.
\end{lemma}

For a vertex $v$ we introduce the projection $\tau_v:B_k\rightarrow A$ to the coordinate $v$.
We denote the $\sigma$-algebra $\tau_v^{-1}(\mathcal{A}(\bA))$ on $B_k$ by $\mathcal{A}_v$.
Obviously $\mathcal{A}_v\subset\mathcal{A}(B_k)$.
By abusing the notation we also define the maps $\tau_v$ on $\bA^{S(0)}$ by
$\tau_v(x)=\tau_v(\delta_k(x))$.
The following lemma shows the connection between the Gowers norm and the cubic structure.

\begin{lemma}\label{gowcub} Let $f:\bA\mapsto\mathbb{C}$ be a bounded measurable function. Then
$$U^{2^k}_k(f)=\int_{x\in B_k}\prod_{v\in V_k}f^{\epsilon(v)}(\tau_v(x))$$
where $f^{\epsilon(v)}$ denotes $\bar{f}$ whenever $|v|$ is odd and denotes $f$ whenever $|v|$ is even.
\end{lemma}

Let $f=(f_v)_{v\in V_k}$ be a collection of bounded measurable functions on $\bA$.
We will need the $2^k$-form
$$\tilde{U}_k(f)=\int_{B_k}\prod_{v\in V_k}f_v(\tau_v(x)).$$
This form naturally extends to the $2^k$-th tensor power of $L_\infty(\bA)$ as a linear function. By abusing the notation we denote that linear function in the same way.

The next lemma shows the connection of the cubic structure with the Fourier $\sigma$-algebras.

\begin{lemma}\label{metszet} For every $v\in V_k$
$$\mathcal{A}_v\cap\langle\{\mathcal{A}_w|w\neq v\}\rangle=\tau_v^{-1}(\mathcal{F}_{k-1}).$$
\end{lemma}

\begin{proof} Let $\mathcal{B}_v$ denote the $\sigma$-algebra $\langle\{\mathcal{A}_w|w\neq v\}\rangle$ and let $\mathcal{B}\subseteq\mathcal{A}(\bA)$ denote the unique $\sigma$-algebra with the property that $\tau^{-1}_v(\mathcal{B})=\mathcal{B}_v$. The transitivity of the symmetry group of the cube on $V_k$ implies that $\mathcal{B}$ does not depend on $v$.
Let $f:\bA\mapsto\mathbb{C}$ be any bounded measurable function.
From lemma \ref{gowcub} we get that
$$U^{2^k}_k(f)=\int_{x\in B_k}\prod_{v\in V_k}E(f^{\epsilon(v)}(\tau_v(x))|\mathcal{B}_v)=U_k(E(f|\mathcal{B})).$$
This implies that $\mathcal{F}_{k-1}\subseteq\mathcal{B}$.

Now we prove that $\mathcal{B}\subseteq\mathcal{F}_{k-1}$. Using the bijection $\delta_k$ it suffices to show that if $f:\bA\mapsto\mathbb{C}$ is a bounded measurable function such that on $\bA^{S(0)}$ the function
$\tau_{[k]}(f)$ is measurable in the $\sigma$-algebra generated by $\{\tau_v|v\in V_k~,~v\neq [k]\}$ then $f\in\mathcal{F}_{k-1}$. This follows easily from theorem \ref{fouequiv}.
\end{proof}

\subsection{Face actions}

Let $\mathcal{F}_{k-1}\subset\mathcal{H}\subset\mathcal{F}_k$ be an $A$ invariant $\sigma$-algebra and let $\sigma\in Aut_k(\mathcal{H})$. Let furthermore $T$ be a subset of $V_{k+1}$. In this subsection we introduce the action $l(T,\sigma)$ on the tensor product
$\otimes^{V_k}L_\infty(\mathcal{H})$ defined by

\begin{equation}\label{edgeact}l(T,\sigma)\bigotimes_{v\in V_{k+1}}f_v=\bigotimes_{v\in V_{k+1}}\sigma_vf_v
\end{equation}

where $\{f_v\}_{v\in V_{k+1}}$ is a collection of bounded $\mathcal{H}$-measurable functions and
$\sigma_v=\sigma$ whenever $v\in T$ and is the identity elsewhere.
The special case where $T$ is a face is of special importance. Such actions will be called {\bf face actions}.

\begin{lemma}If $e$ is an edge of $V_{k+1}$ and $\sigma\in Aut_k(\bA)$ then action $l(e,\sigma)$ preserves the form $\tilde{U}_{k+1}$.
\end{lemma}

\begin{proof}Without loss of generality we can assume that $e=\{\emptyset,\{k+1\}\}$. Let $\{f_v\}_{v\in V_{k+1}}$ be a collection of bounded real valued functions on $\bA$.. Using lemma \ref{spidrep}we obtain that
$$\tilde{U}_{k+1}(l(e,\sigma)\otimes_{v\in V_{k+1}}f_v))=$$
$$=\int_{a_0,a_1,\dots,a_{k+1}}\sigma f_\emptyset(a_0)\sigma f_{\{k+1\}}(a_0+a_{k+1})\prod_{v\in V_{k+1}\setminus e}f_v(a_0+\sum_{i\in v} a_i)=$$
$$=\int_{a_{k+1}}\int_{a_0,a_2,\dots,a_{k+1}}\sigma(f_\emptyset(a_0)f_{\{k+1\}}(a_0+a_{k+1}))\prod_{v\in V_{k+1}\setminus e}f_v(a_0+\sum_{i\in v} a_i)$$
The inner integral can be written as
\begin{equation}\label{vetit}
\tilde{U}_k(\sigma \hat{f}_\emptyset\otimes\bigotimes_{v\in V_k\setminus\emptyset}\hat{f}_v)
\end{equation}
where
$$\hat{f}_v(x)=f(x+a_{k+1})f(x).$$
Now lemma \ref{metszet} implies that in (\ref{vetit}) the term $\sigma\hat{f}_\emptyset$ can be replaced by
$$E(\sigma\hat{f}_\emptyset|\mathcal{F}_{k-1})=E(\hat{f}_\emptyset|\mathcal{F}_{k-1})$$
and then we can replace it by $\hat{f}_\emptyset$.
After this transformation (\ref{vetit}) becomes
$$\tilde{U}_k(\bigotimes_{v\in V_k}\hat{f}_v).$$
This finishes the proof.
\end{proof}

\subsection{Isometric extensions are Abelian}

The main result in this section is the following lemma.

\begin{lemma}[Abelian extension]\label{abext} Let $f:\bA\rightarrow G$ be a minimal pre-cocycle of order $k$. Then $G$ is Abelian.
\end{lemma}

\begin{proof} Let $e_1$ and $e_2$ be two edges of $V_{k+1}$ intersecting each other at a vertex $w\in V_{k+1}$.
Let $g_1,g_2$ be two elements from $G$. Let $\mathcal{H}$ be the $\sigma$ algebra guaranteed by lemma \ref{cosigma}. Now we denote by $c_{w}$ the commutator $[l(e_1,\varrho(g_1),l(e_2,\varrho(g_2)]$. We have that $c_{w}$ is acting on $\otimes^{V_{k+1}}L_\infty(\mathcal{H})$ and the action is given by
$$c_w=l(w,\varrho([g_1,g_2])).$$
From lemma... we ga that the form $\tilde{U}_{k+1}$ is preserved by $c_w$.
This means that the composition of such actions for a sequence of different vertices $w$ also preserves $\tilde{U}_{k+1}$.
From this we deduce that $\varrho([g_1,g_2])$ has to be trivial.
Let $f$ be a bounded real valued $\mathcal{H}$ measurable function and $g=f-\varrho([g_1,g_2])f$. Then
$$U_{k+1}(g)=\int_x\prod_{v\in V_{k+1}}(f(\tau_v(x))-\varrho([g_1,g_2])f(\tau_v(x))).$$
Expanding the product on the right side all the terms have the same integral and the number of positive signs is the same as the number of negative signs. We get that $U_{k+1}(g)=0$.
Since $g\in\mathcal{F}_k$ this implies that $g=0$.
\end{proof}

\section{Higher order Fourier analysis}

\subsection{Definition of higher order Fourier analysis}

The $i$-th order Fourier analysis deals with the decomposition of the module $L_2(\mathcal{F}_k)$ into irreducible $\bA$ invariant $L_\infty(\mathcal{F}_{k-1})$-modules. The fundamental theorem of higher order fourier analysis is the following.

\begin{theorem}[Fundamental Theorem] The rank one $L_\infty(\mathcal{F}_{k-1})$ modules are pairwise orthogonal and they generate the Hilbert space $L_2(\mathcal{F}_k)$. Every such rank one module has a generator of the form
$\alpha:\bA\rightarrow\mathbb{C}$ with $|\alpha|=1$ and $\alpha(x)\overline{\alpha(x+t)}\in\mathcal{F}_{k-1}$. These functions are called $k$-th {\bf order characters}.
\end{theorem}

It will turn out that the above mentioned rank $1$ modules are forming an Abelian group using the point-wise multiplication.

\subsection{Higher order group algebras}

In this section we introduce algebras which replace convolution algebras on ordinary compact topological groups in the higher order setting.
Let first $M_k$ denote the space $L_2(\bA\times\bA,\mathcal{F}_k\otimes\mathcal{F}_k,\bm\times\bm)$ of Hilbert-Schmidt operators.
It is easy to check that $M_k$ has a Banach $C^*$-algebra structure with multiplication
$$(K_1\circ K_2)(x,y)=\int_z K_1(x,z)K_2(z,y)~d\bm.$$

\begin{definition}[Higher order group algebras] Let $C_k\subset M_k$ denote the set of elements $K\in M_k$ such that the functions
$$f_{x,y}(t):=K(x+t,y+t)$$
are in $L_2(\mathcal{F}_{k-1})$ for every $x,y\in\bA$.
\end{definition}

\begin{lemma} The set $C_k$ is forming a sub $C^*$-algebra of $M_k$.
\end{lemma}

\begin{proof} Use lemma \ref{inth}.
\end{proof}

Now we describe a useful way of constructing elements in $D_k$.

\begin{lemma} Let $f,g$ be two functions in $L_2(\mathcal{F}_k)$. Then there is an element $e=e(f,g)$ in $C_k$ such that for every fixed $x,y$
$$e(x+t,y+t)=E(f(x+t)g(y+t)|\mathcal{F}_{k-1})(t).$$
\end{lemma}

\subsection{Decomposition of the elements in $C_k$}

Let $C$ be a self adjoint element of $C_k\cap L_\infty(\bA\times\bA)$. We denote by $im(C)\subset L_2(\bA)$ the image space of the operator $C$.
Note that $im(C)$ is a separable Hilbert subspace of the non separable space $L_2(\bA)$.
It is classical that
$$im(C)=\bigoplus V_i$$
where $V_i$ are finite dimensional eigenspaces of $C$ corresponding to different eigenvectors.
The spaces $V_i$ are contained in $L_\infty(\bA)$.
The main result of this section is the following.

\begin{lemma}\label{vegesdim} Each space $V_i$ is contained in a finite rank module over $L\infty(\mathcal{F}_{k-1})$. Furthermore this finite rank module can be decomposed into rank $1$ modules.
\end{lemma}

\begin{proof}
Let $P_i$ denote the orthogonal projection to $V_i$. The operators $P_i$ are contained in the Banach algebra generated by $C$ so they are all elements of $C_k$.
Let $f_1,f_2,\dots,f_d\in L_\infty(\bA)$ be an orthonormal basis for $V_i$.
We also introduce the vector valued function $f:\bA\mapsto \mathbb{C}^d$ defined by
$f(x)=(f_1(x),f_2(x),\dots,f_d(x))$.
The operator kernel $P_i$ satisfies $P_i(x,y)=\langle f(x),f(y)\rangle$.
Since $P_i\in C_k$ we get that
$\langle f(x),f(x+t)\rangle\in \mathcal{F}_{k-1}$ for every fixed $t\in\bA$.

We show that there are elements $t_1,t_2,\dots,t_d$ such that
the function
$\rk(f(x+t_1),f(x+t_2),\dots,f(x+t_d))$ is equal to $d$ on a positive measure set in $\mathcal{F}_k$.
Let $$h(t_1,t_2,\dots,t_d):=|det(f(t_1),f(t_2),\dots,f(t_d))|^2.$$ It is easy to see that the linear independence of $f_1,f_2,\dots,f_d$ implies that $h>0$ on a positive measure subset set in $\bA^d$. This implies that there is a fixed vector $t_1,t_2,\dots,t_d$ such that $h(t_1+x,t_2+x,\dots,t_d+x)>0$ on a positive measure set of $\bA$.
Since the matrix $W_{i,j}(x)=\langle f(t_i+x),t(t_j,x)\rangle$ is measurable in $\mathcal{F}_{k-1}$ we get that the set where $$h(t_1+x,t_2+x,\dots,t_d+x)=det(W_{i,j}(x))$$ is positive is on $P\in\mathcal{F}_{k-1}$.

The next step is that we run the Gram-Schmidt orthogonalization on the set $P$ for the vectors $\{f(x+t_i)\}_{i=1}^d$.
It is easy to see that all the coefficients are measurable in $\mathcal{F}_{k-1}$ so we find functions
$o_1,o_2,\dots,o_d:P\rightarrow \mathcal{C}^d$ such that
$\{o_i(x)\}_{i=1}^d$ is an orthonormal basis for every $x\in P$ and furthermore
$$o_i=\sum_{j=1}^d \lambda_{i,j}(x)f(x+t_j)$$ for every $1\leq i\leq d$ where $\lambda_{i,j}(x)$ is measurable in $\mathcal{f}_{k-1}$.
From this we get that every translate $f(x+t)$ of $f$ can be expressed on $P$ from $\{o_i\}_{i=1}^d$ using $\mathcal{F}_{k-1}$ measurable functions as coefficients.

Now we use Rohlin's lemma to find at most countable many subsets $P_1,P_2,\dots$ of $P$ all measurable in $\mathcal{F}_{k-1}$ such that $\bA$ is a disjoint union of translated versions of these sets.
By translating the basis $\{o_i\}_{i=1}^\infty$ correspondingly one obtains an extension of the orthonormal basis $\{o_i\}_{i=1}^\infty$ to the whole set $\bA$.
In this extended basis any translate of $f$ can be expressed with $\mathcal{F}_{k-1}$ measurable coefficients.
This shows that $V_i$ is contained in a module of rank at most $d$.

Let $O(x)$ denote the unitary matrix formed by the $\{o_i(x)\}_{i=1}^d$. It is clear that $x\rightarrow O(x)$ is a pre cocycle of degree $k$. This means by lemma \ref{findmin} and lemma \ref{abext} that it is equivalent with a cocycle $O'(x)$ going to some abelian subgroup $H$ of the unitary group. By decomposing $H$ into $1$ dimensional irreducible representations we find pre-cocycles $O_i'(x)$ taking values in $\mathbb{C}$ of absolute value $1$.
This provides rank one modules over $\mathcal{F}_{k-1}$ that generate $V_i$.

\end{proof}

\begin{corollary}\label{rkone} Let $C\in C_k\cap L_\infty(\bA\times \bA)$ be an arbitrary element. Then $im(C)$ is contained in a space generated by rank $1$ modules.
\end{corollary}

\begin{proof} Using that $im(C^*C)=im(C)$ we can assume that $C$ is self adjoint. Then the previous lemma completes the proof.
\end{proof}

\subsection{Proof of the fundamental theorem}

Let $\mathcal{F}_{k-1}\subset \mathcal{H}\subset \mathcal{F}_k$ be an $\bA$ invariant $\sigma$ algebra that is relative separable with respect to $\mathcal{F}_{k-1}$.
Let furthermore $f_1,f_2,\dots$ be a relative orthonormal basis of $L_2(\mathcal{H})$ over $\mathcal{F}_{k-1}$. This means that $E(f_if_j|\mathcal{F}_{k-1})$ is the constant $1$ function for $i=j$ and is the constant $0$ function for $i\neq j$.
Every element $f\in L_2(\mathcal{H})$ can be uniquely written as
$\sum_i \lambda_i(x)f_i$ where the functions $\lambda_i$ are $\mathcal{F}_{k-1}$ measurable and can be obtained by
$\lambda_i=E(ff_i|\mathcal{F}_{k-1})$.
Now let us introduce the uniquely determined functions $\lambda_{i,j}(t,x)$ by
\begin{equation}\label{eltolas}
f_i(x+t)=\sum_j \lambda_{i,j}(t,x)f_j(x).
\end{equation}
Let furthermore $M_{x,y}$ denote the $\mathbb{N}\times\mathbb{N}$ matrix defined by
$$M_{x,y}(i,j)=\lambda_{i,j}(y-x,x).$$
Lemma \ref{projprop} guarantees that $M_{x,y}$ is a measurable function of $(x,y)$.
It is easy to check from (\ref{eltolas}) that
$$M_{x,y}M_{y,z}=M_{x,z}$$
hold for almost all triples $x,y,z$
$$M_{x,y}=M_{y,x}^{-1}=M_{y,x}^*$$
holds for almost all pairs $x,y$.
This implies that there is a fixed $y$ such that the above equations hold for almost all $x$ and $z$.
Let $M_x:= M_{x,y}$. We have that $M_{x,z}=M_xM_z^{-1}$ for almost all $x,z$.
By definition we have that $M_xM_{x+t}^{-1}$ is measurable in $\mathcal{F}_{k-1}$ for every fixed $t$.
This guarantees that $M_k$ is measurable in $\mathcal{F}_k$ since the $\sigma$-algebra generated by $M_x$ has to be relative separable with respect to $\mathcal{F}_{k-1}$.
Let $\mathcal{H}'$ be the $\sigma$-algebra generated by the function $x\rightarrow M_x$.
It is easy to see that $\mathcal{H}\subseteq\mathcal{H}'$.
Summarizing the above information we obtain that $M_{x,y}(i,j)\in C_k$ for every fixed $i,j$. Using Corollary \ref{rkone} we get that $M_{x,y}$ is measurable in a $\sigma$ algebra generated by rank $1$ modules.
So the $\sigma$-algebra $\mathcal{H}'$ is generated by rank $1$ modules.
Since $\mathcal{F}_k$ is the union of all relative $\mathcal{F}_{k-1}$ separable invariant $\sigma$-algebras the proof is complete.

\vskip 0.2in

\noindent
Bal\'azs Szegedy
\noindent
University of Toronto, Department of Mathematics,
\noindent
St George St. 40, Toronto, ON, M5R 2E4, Canada

\end{document}